\documentclass[a4paper,11pt]{article}
\usepackage{amsmath,amssymb,xspace}
\usepackage{mathtools}
\date {}
\def \Z {\mathbb{Z}}
\def \C {\mathbb{C}}
\def \F {\mathbb{F}}
\newcommand{\ab}{\mathrm{ab}}
\newcommand{\Hom}{\mathrm{Hom}}
\newcommand{\SL}{\mathrm{SL}}
\newcommand{\GL}{\mathrm{GL}}
\newcommand{\Syl}{\mathrm{Syl}}
\newcommand{\Aut}{\mathrm{Aut}}

\newtheorem{theorem}{Theorem}[section]
\newtheorem{proposition}[theorem]{Proposition}

\usepackage{latexsym}
\title{Finite abelian groups are Schur multipliers of finite solvable groups}
\author{Derek F.\ Holt\\
Mathematics Institute, University of Warwick,\\
Coventry CV4 7AL, UK. D.F.Holt@warwick.ac.uk}
\begin{document}
\maketitle
\begin{abstract} We prove that every finite abelian group arises as the
Schur multiplier of a finite solvable group, from which it can be deduced that
every finite abelian group is the knot group of some extension of number
fields.
\end{abstract}
\section{Introduction}
\label{sec:intro}
We prove the following result.
\begin{theorem}
Every finite abelian group arises as the
Schur multiplier of a finite solvable group.

More precisely, for every finite abelian group $A$, there is a finite solvable
group $G$ for which the Schur multiplier $M(G)$ of $G$ is isomorphic to $A$.
Furthermore we can choose such a group $G$ with a normal nilpotent subgroup
$N$ of class at most $2$ such that $G/N$ is cyclic.
\end{theorem}

Some motivation for this result comes from the study of the \emph{knot group}
$\mathfrak{K}(K/k)$ of an extension $K/k$ of number fields. This was
originally defined by Scholz in~\cite{Scholz36}, and it is called the
\emph{number knot} in~\cite{Jehne79}.
For $c \in k^\times$ (the multiplicative group $k \setminus \{0\}$), we define
$X_c$ to be the norm equation $N_{K/k}(x_1w_1 + \cdots x_dw_d) = c$, where
$\{w_1,\ldots,w_d\}$ is a $k$-basis of $K$. Then the knot group
$\mathfrak{K}(K/k)$ is isomorphic to be the quotient $X/Y$, where $X$ and $Y$
are respectively the subgroups of $k^\times$  consisting of those elements $c$
for which the equation $X_c$ is solvable locally for all primes of $k$, and
those for which $X_c$ is solvable globally. It can be deduced from a
result of Shafarevich~\cite[Proposition 5.5]{FLN18} that, if $M(G)=A$ for some
finite solvable group $G$, then $A$ arises as a knot group of some
extension $K/k$ of number fields. So the theorem above implies that all finite
abelian groups arise as knot groups.

We start by summarizing the properties of the Schur multiplier $M(G)$ of a
(finite) group $G$ that we shall be needing in the proof of the theorem.
 
$M(G)$ is usually defined to be the second homology group $H_2(G,\Z)$ of $G$
with trivial action on $\Z$. Alternatively,
for finite groups $G$, as in for example \cite[V:\S 23]{Huppert67},
$M(G)$ is sometimes defined to be $H^2(G,K^\times)$, where $K^\times$ is the
multiplicative group of an algebraically closed field of characteristic $0$,
such as $\C$. By the Universal Coefficient Theorem, this is equal to the
dual $\Hom(M(G),K^\times)$ of $M(G)$, which is isomorphic to $M(G)$ when it
is finite. By \cite[V: 23.5]{Huppert67}, $M(G)$ is finite when $G$ is finite.

A group $\widehat{G}$ with normal subgroup $Z \le Z(G) \cap [G,G]$ is said
to be a \emph{stem extension} of $G$. Any stem extension is an epimorphic
image of a maximal stem extension, also known as a \emph{Schur cover} of $G$.
For a Schur cover $\widehat{G}$ of $G$, we have $Z \cong M(G)$; see, for
example, \cite[V: 23.5]{Huppert67}

The following theory can be found in~\cite[Chap XII, Sec 10]{CE56}.
For an abelian group $A$ and a prime $p$, we denote the $p$-primary part of
$A$ by $A_p$.  If $G$ is finite, $P \in \Syl_p(G)$, and $M$ is a finitely
generated
$\Z G$-module, then the cohomology group $H^i(G,M)$ is finite for each $i>0$,
and the restriction map $H^i(G,M)_p \rightarrow H^i(P,M)$ is a monomorphism. It
is proved in~\cite[Chap XII, Thm 10.1]{CE56} that the image of this map is equal
to the subgroup of $H^i(P,M)$ consiting of those elements that are
\emph{stable} with respect to $G$.
We shall not need the precise definition of stable elements here; we need
only the case when $P \unlhd G$, in which case the stable elements
consist of those elements of $H^i(P,M)$ that are fixed under the natural
induced action of $N_G(P)$ on $H^i(P,M)$.

In the case when $M = \C^\times$ with trivial action of $G$ on $M$,
as we saw above, we have $H^2(G,\C^\times) \cong \Hom(M(G),\C^\times)$,
the dual of $M(G)$. In this situation, if $P \unlhd G$, then
$M(G)_p$ is isomorphic to the $M(P)/[M(P),G]$ with the natural induced action of
$G$ on $M(P)$, and we shall use this formula to compute $M(G)_p$.

For convenience we shall state some further standard results on the Schur
multiplier. See, for example, \cite[V:25.3, 25.10]{Huppert67} for proofs.
Since finite abelian groups are direct products of cyclic groups, the
third part follows from the first two. For the final statement, note
that the results in~\cite{CE56} mentioned earlier imply that
restriction map $H^2(G,\C^\times)_p \to H^2(H,\C^\times)_p$
is injective, and that $H^2(G,\C^\times)_p$ and $H^2(H,\C^\times)_p$ are
isomorphic to the duals of $M(G)_p$ and $M(H)_p$.
We denote the abelianization $G/[G,G]$ of a group $G$ by $G_\ab$.

\begin{proposition}\label{prop:dp}
\begin{enumerate}
\item[(i)] For a cyclic group $G$ we have $M(G)=1$.
\item[(ii)] For a direct product $G \times H$ of groups $G$ and $H$, we have
$M(G \times H) \cong M(G) \times M(H) \times G_\ab \otimes H_\ab$.
\item[(iii)] For a finite abelian group $A$, we have $M(A) = A \wedge A$, the
\emph{wedge product} of $A$ with itself.
\item[(iv)] If $H \le G$ with $G$ finite and $|G:H|$ is not divisible by the
prime $p$ (in particular when $P \in \Syl_p(G)$), then $M(G)_p$ is isomorphic
to a quotient group of $M(H)_p$.
\end{enumerate}
\end{proposition}

We shall also need the following result, known as \emph{Fitting's theorem},
proved in \cite[Thm 4.34]{Isaacs08}.
\begin{proposition}\label{autpgp}
It $T$ is a $p'$-group of automorphisms of a finite abelian $p$-group $A$,
then $A = C_A(T) \times [A,T]$.
\end{proposition}

We turn now to the proof of the theorem.  If $A$ is trivial then we can choose
$G$ to be trivial, so we may assume that $A$ is non-trivial.
In the following three sections, we prove the theorem in the cases when $A$ is
an elementary abelian $p$-group, an arbitrary finite abelian $p$-group, and an
arbitrary finite abelian group.

\section{The elementary abelian case}\label{sec:elab}
Suppose in this section that $A$ is an elementary abelian $p$-group of order
$p^k$ for some prime $p$ with $k \ge 1$.

Let $N = N_1 \times N_2$ with $N_1 \cong N_2 \cong A$, and let
$G = N \rtimes \langle t \rangle$, where $t$ has order $p^k-1$ and acts on $N$
as follows. We identify $N_1$ and $N_2$ with the additive groups of the
field $\F_{p^k}$ of order $p^k$, let $x$ be a generator of $\F_{p^k}^\times$,
and let $t$ act on $N_1$ and $N_2$ by multiplication by $x$ and $x^{-1}$
respectively.

We claim that, except in the case when $p^k = 4$, we have $M(G) \cong A$.
Note that, since the Sylow $q$-subgroups $Q$ of $G$ are cyclic for primes
$q \ne p$, it follows from parts (i) and (iv) of Proposition~\ref{prop:dp}
that $M(G)_q=1$, so $M(G)=M(G)_p = M(N)/[M(N),G] =
M(N)/[M(N),\langle t \rangle]$, which we shall denote by $M(N)/[M(N),t]$. 

By Proposition~\ref{prop:dp} \,(iii),
we have $M(N) \cong M(N_1) \times M(N_2) \times N_1 \otimes N_2$.
Since $t$ fixes $N_1$ and $N_2$, it fixes the subgroups
$M(N_1)$, $M(N_2)$, and $N_1 \otimes N_2$ in its induced
action on $M(N)$, so we can consider these three induced actions separately.
The following calculations are similar to those in the proof of
\cite[V:25.5]{Huppert67} (which states that the $p$-part of $M(\SL(2,q))$
is trivial when $q \ne 4,9$).

We can regard $N$ and its direct factors as vector spaces over $\F_p$, and
the action of $t$ on $N$ induces a linear map $\tau:N \to N$. In general,
it is easy to see that, if $\tau_1:V_1 \to V_1$ and $\tau_2:V_2 \to V_2$ are
linear maps of vector spaces $V_1$ and $V_2$ over the same field, then the
eigenvalues (in a suitable extension field) of the linear map
$\tau_1 \otimes \tau_2$ of $V_1 \otimes V_2$ are
the products of the eigenvalues of $\tau_1$ and $\tau_2$. Similarly, if the
eigenvalues of the linear map $\tau:V \to V$ are $\lambda_1,\ldots,\lambda_k$,
then those of $\tau$ in its induced action on $V \wedge V$ are
$\{ \lambda_i\lambda_j : 1 \le i < j \le k \}$.

We show first that $[M(N_1),t] = M(N_1)$ when $p^k \ne 4$ which, by
Proposition~\ref{autpgp}, is equivalent to showing that $C_{M(N_1)}(t) = 1$
or, equivalently, that the linear map $\tau:N_1 \to N_1$ induced by the
action of $t$ has no eigenvalues equal to $1$ in its induced
action on $M(N_1)$. Since $\tau$ acts on $N_1$ by multiplication by $x$ in
the field $\F_{p^k}$, we see that the eigenvalues of $\tau$ in $\F_{p^k}$ are
the Galois conjugates of $x$ in $\F_{p^k}$, which are
$x^{p^i}$ for $0 \le i < k$. So the eigenvalues of $\tau$ on $V \wedge V$
are $x^{p^i+p^j}$ for $0 \le i < j < k$. Such an eigenvalue would be equal to
$1$ if and only if $p^i+p^j \equiv 0 \pmod{p^k-1}$. But, except in the
case $p=2$, $i=0$, $j=1$, $k=2$, we have $p^i+p^j < p^k-1$, which proves 
our claim. The proof that $[M(N_2),t] = M(N_2)$ is the same with $x$
replaced by $x^{-1}$. 

We complete the proof of the theorem in this case when $p^k \ne 4$ by proving
that $M(N_1 \otimes N_2)/[N_1 \otimes N_2,t]$ has order $p^k$ or, equivalently
using Proposition~\ref{autpgp}, that $\tau$ has exactly $k$ eigenvalues equal
to $1$ in its action on $N_1 \otimes N_2$. Since this action is the same as
that of $\tau_1 \otimes \tau_2$, where $\tau_1$ acts of $N_1$ by multiplication
by $x$ and $\tau_2$ acts on $N_2$ by multiplication by $x^{-1}$, we see that the
eigenvalues of the induced action of $\tau$ are
$x^{p^i-p^j}$ for $0 \le i,j < k$. So the cases $i=j$ give rise to $k$
eigenvalues equal to $1$. Since we have $|p^i-p^j| < p^k-1$ for all $p$
and $0 \le i,j < k$, these are the only eigenvalues equal to $1$.

For the exceptional cat $p^k=4$, we again choose
$G = N \rtimes \langle t \rangle$ with $|N| = p^{2k} = 16$, but now we choose
$t$ to have order $5$ and act irreducibly on $N$. Then the induced action of
$\tau$ of $t$ on $N \wedge N$ decomposes as a sum of two trivial modules and an
irreducible  module of dimension $4$ (this is easily verified by computer),
so we get the required result.

Note that in all cases, since $t$ acts without fixed points on $N$, we have
$[G,G] = N$ with $G_\ab \cong \langle t \rangle$.

\section{The $p$-group case}\label{sec:pgrp}
Suppose first that $A$ is a homocyclic $p$-group $C^k$, where
$C$ is cyclic of order $p^n$. We dealt with the case $n=1$ in
Section~\ref{sec:elab}, so we can assume that $n>1$.
It is well-known (see, for example, \cite{HR07}) that $\Aut(A)$ has a normal
$p$-subgroup $B$  with $\Aut(A)/B \cong \GL(k,p)$, where $B$ induces the
trivial action on $A/A^p$, and $\Aut(A)/B$ induces the
natural action of $\GL(k,p)$ on $A/A^p$.

In particular, there is an automorphism $\alpha$ of order $p^k-1$ of $A$ that
induces the automorphism of $A/A^p$ considered in Section~\ref{sec:elab},
which acts by multiplication by a generator $x$ of $\F_{p^k}$.

Let $N = N_1 \times N_2$ with $N_1$ and $N_2$ homocyclic and isomorphic to
$A$, and let $G = N \rtimes \langle t \rangle$ where $t$ acts as $\alpha$
on $N_1$ and as $\alpha^{-1}$ on $N_2$. We have
$M(N) = M(N_1) \times M(N_1) \times N_1 \otimes N_2$, where $M(N_1)$
and $M(N_2)$ are isomorphic to $A\wedge A$, which is homocyclic and
isomorphic to $C^{k(k-1)/2}$, and $N_1 \otimes N_2$ is homocyclic and
isomorphic to $C^{n^2}$.
Again we have $M(G)=M(G)_p = M(N)/[M(N),t]$.

The induced action of $t$ on $M(N)/M(N)^p$ is the
same as its action on $M(N)$ in Section~\ref{sec:elab} when $N$ was elementary
abelian and so, when $p^k \ne 4$, the subgroup of $M(N)/M(N)^p$ fixed
by $t$ in this action has order $p^k$.
From Proposition~\ref{autpgp} we have $M(N) = C_{M(N)}(t) \times [M(N),t]$
which, together with the known induced action of $t$ on $M(N)/M(N)^p$,
implies that $C_{M(N)}(t)$ and $[M(N),t]$ are homocyclic with
$C_{M(N)}(t)$ isomorphic to $C^k$, which gives the desired result in this case.
When $p^k=4$, we use the same variant as in Section~\ref{sec:elab} and choose
$N= C^4$ with $t$ of order $5$.

As we noted at the end of Section~\ref{sec:elab}, $t$ acts without fixed points
on $N/N^p$ and so we see from Proposition~\ref{autpgp} that $[t,N]=N$,
so $[G,G] = N$ in all cases.

Now suppose that $A$ is an arbitrary abelian $p$-group of exponent $p^n$ with
$n>1$ and let $k$ be the rank of $A$.
Let $H = N \rtimes \langle t \rangle$ with $M(H) \cong C^k$
and $C$ cyclic of order $p^n$ be defined as in the case when $A$ is
homocyclic, and let $\widehat{H}$ be a Schur cover of $H$.
So there exists $K \le Z(\widehat{H}) \cap [\widehat{H},\widehat{H}]$
with $K \cong M(H) \cong C^k$ and $\widehat{H}/K \cong H$.
Let $L$ be the complete inverse image of $N$ in $\widehat{H}$. Then $L$ is a
$p$-group of class $2$, and we have $\widehat{H} = L \rtimes \langle t \rangle$,
where we are abusing notation by using $t$ to denote both an element
in $H$ and an inverse image in $\widehat{H}$ having the same order.

Now $K$ has a subgroup isomorphic to $A$, which we shall also denote by $A$.
We claim that $\widehat{H}$ is a Schur cover of $\widehat{H}/A$, which implies
that $A \cong M(\widehat{H}/A)$. So this will prove the theorem with
$G = \widehat{H}/A$, which is a semidirect product of the $p$-group 
$L/A$ with $\langle t \rangle$.

If not, then $\widehat{H}$ is a proper quotient of a Schur cover $F$ of
$\widehat{H}/A$, where $F$ has non-trivial
central subgroups $Z < \widetilde{A}$ with $F/Z \cong \widehat{H}$,
and $\widetilde{A} \cong M(\widehat{H}/A)$ the complete inverse image
of $A$ in  $F$. As before, all Sylow $q$-subgroups of $\widehat{H}/A$ for
$q \ne p$ are cyclic, so $M(\widehat{H}/A)$ and hence also $Z$ are $p$-groups.

Let $\widetilde{L}$ and $\widetilde{K}$ be the complete inverse images of $L$
and $K$ in $F$, and note that $Z < [F,F]$ and $K < [\widehat{H},\widehat{H}]$
together imply that $\widetilde{K} < [F,F] =  \widetilde{L}$.
Furthermore, we have $F = \widetilde{L} \rtimes \langle t \rangle$,
where we are again abusing notation and re-using $t$ as an element of the
same order as the original $t$.

We shall show that $\widetilde{K} \le Z(F)$ which, together with
$\widetilde{K} \le [F,F]$ and  $F/\widetilde{K} \cong \widehat{H}/K
\cong H$, implies that $F$ is a larger stem extension of $H$ than $\widehat{H}$,
a contradiction. 

We have $[[F,\widetilde{K}],F] \le [Z,F]=1$ and so, by the three subgroups
lemma, we get $[[F,F],\widetilde{K}] = [\widetilde{L}, \widetilde{K}] = 1$.
But the $p'$-element $t$ induces the identity automorphism on the $p$-groups
$\widetilde{K}/Z$ and $Z$, so it centralizes $\widetilde{K}$ (see, for example,
\cite[Lemma 4.29]{Isaacs08}).
So, since $F = \langle \widetilde{L}, t \rangle$, we have
$\widetilde{K} \le Z(F)$, as claimed. This completes the proof of the
theorem when $A$ is a $p$-group.

Now in all cases we have defined a group $G$ with $M(G) = A$ and, after
renaming $L/A$ as $N$, we have $G = N \rtimes \langle t \rangle$ with
$N=[G,G]$ nilpotent of class at most $2$, and $t$ of order $p^k-1$
(or $5$ when $p^k=4$). 

We observe that, for any $m>0$ that is coprime to $p$, we can replace
$t$ by an element of order $tm$ that induces the same action as the original
$t$ on $N$, where $t^{p^k-1}$ (or $t^5$ when $p^k=4$) centralizes $N$,
and we will still have $M(N \rtimes \langle t \rangle) \cong A$, with the
proof unchanged.

\section{The general case}\label{sec:general}
Now assume that $A$ is a direct product of $p$-groups $A_{p_i}$ of order
$p_i^{k_i}$ for distinct primes $p_i$, and for each $i$ let
$G_i = N_i \rtimes \langle t_i \rangle$ be the group defined in the previous
sections with $M(G_i) = A_{p_i}$.

We define $G := N \rtimes \langle t \rangle$, where $N$ is the direct
product of the groups $N_i$ and $t$ is an element of order the least
common multiple of the orders $|t_i|$ of the elements $t_i$.
(So $|t_i| = p_i^{k_i} - 1$, or $5$ when $p_i^{k_i}=4$.)
The action of $t$ on each $N_i$ is defined to be the induced action of
$t_i$ with $t^{|t_i|}$ centralizing $N_i$.

We complete the proof of the theorem by showing that $M(G) \cong A$. Since
the Sylow $q$-subgroups of $G$ are cyclic for primes $q$ unequal to any $p_i$,
we have $M(G)_q=1$ for such primes, and it is enough to show that
$M(G)_{p_i} = A_{p_i}$ for each of the primes $p_i$. So let $p=p_i$
for some fixed $i$.

Let $N'$ be the direct product of the groups $N_j$ for $j \ne i$. Then $G$ is
the semidirect product $N' \rtimes H$, where
$H = N_i \rtimes \langle t \rangle$.
Factor the order $|t|$ of $t$ as $ab$, where the prime factors of $a$ all
divide $|t_i|$, and those of $b$ do not divide $|t_i|$. Then we can write
$t = rs$, with $|r|=a$ and $|s|=b$ and, since $s$ centralizes $N_i$,
the group $H$ is the direct product
$(N_i \rtimes \langle r \rangle) \times \langle s \rangle$.
Then, since the abelianizations of the direct factors have coprime orders $a$
and $b$, Proposition~\ref{prop:dp}\,(i),\,(ii) gives
$M(H) \cong M(N_i \rtimes \langle r \rangle)$.
Since $a = |r|$ is coprime to $p$, and $a$ is divisible by $|t_i|$, it
follows from the remark at the end of the previous section that
$M(H) \cong M(G_i) \cong A_{p_i}$.   

Since $|G:H| = |N'|$ is coprime to $p$, $M(G)_p$ is isomorphic to a quotient
group of $M(H)$ by Proposition~\ref{prop:dp}\,(iv). But a Schur cover of $H$
lifts to a stem extension of $G$, so $|M(G)_p| \ge |M(H)|$ and hence
$M(G)_p \cong M(H) \cong A_{p_i}$ as required.
This completes the proof of the theorem.

\bibliographystyle{alpha}
\bibliography{solmult}
\end{document}